\documentclass[letterpaper,conference]{IEEEtran}
\IEEEoverridecommandlockouts

\usepackage{cite}
\usepackage{amsmath,amsthm,amssymb,amsfonts}
\usepackage{algorithmic}
\usepackage{graphicx}
\usepackage{textcomp}
\usepackage{xcolor}
\def\BibTeX{{\rm B\kern-.05em{\sc i\kern-.025em b}\kern-.08em
    T\kern-.1667em\lower.7ex\hbox{E}\kern-.125emX}}

\usepackage{dsfont}

\usepackage{verbatim}
\usepackage[T1]{fontenc}
\usepackage[latin9]{inputenc}
\usepackage{graphicx}
\usepackage[english]{babel}
\usepackage{placeins}
\usepackage{comment}
\usepackage{booktabs}
\usepackage{url}
\usepackage{hyperref}
\usepackage{tikz}
\usepackage{mathtools}

\usepackage{babel}

\usepackage{color}
\renewcommand{\S}[0]{Section}

\newtheorem{thm}{Theorem}

\begin{document}

\title{Three Formulations of the Kuramoto Model \\ as a System of Polynomial Equations\thanks{D.M. was supported by NSF-ECCS award ID 1509036. J. M. received funding from the European Union Horizon 2020 Programme (Horizon2020/2014-2020), under grant agreement number 688380. M.N. was supported in part by the Fields Institute for Research in Mathematical Sciences.
The authors can be contacted at:
\url{tchen1@aum.edu}, 
\url{jakub.marecek@ie.ibm.com}, 
\url{dhagashbmehta@gmail.com}, and
\url{matthew.niemerg@alephzero.org}, respectively. 
}} 

\author{\IEEEauthorblockN{
Tianran Chen} 
\IEEEauthorblockA{Auburn University \\ Montgomery, AL, USA} \and
\IEEEauthorblockN{Jakub Mare{\v c}ek} 
\IEEEauthorblockA{IBM Research -- Ireland \\ Dublin, Ireland}
\and
\IEEEauthorblockN{Dhagash Mehta}
\IEEEauthorblockA{University of Notre Dame \\ South Bend, IN, USA}
\and
\IEEEauthorblockN{Matthew Niemerg}
\IEEEauthorblockA{IBM Center of Excellence \\
and Aleph Zero}
}

\maketitle

\pagestyle{headings}

\begin{abstract} 

We compare three formulations of stationary equations of the Kuramoto model as systems of polynomial equations. In the comparison, we present bounds on the numbers of real equilibria
based on the work of Bernstein, Kushnirenko, and Khovanskii,
and performance of methods for the optimisation over the set of equilibria based on the work of Lasserre, both of which could be of independent interest.

 
\end{abstract}



\section{Introduction}\label{Sec:Introduction}
The Kuramoto model is a prototypical model for studying many phenomena including the synchronization  
of power systems, neural networks, chemical oscillators, particle coordinations, 
rhythmic applause, and so on.
See \cite{kuramoto1975self,acebron2005kuramoto,strogatz2000kuramoto,FD-FB:13b}.
The general Kuramoto model is given by the differential equation:
\begin{equation}\label{Eq:KuramotoStandard}
    \frac{d\theta_{i}}{dt} =\omega_{i} - \frac{1}{N} 
    \sum_{j=1}^{N}K_{i,j}
    \sin(\theta_{i}-\theta_{j}), \text{ for } i=1,...,N,
\end{equation}
where $N$ is the number of oscillators, $K_{i,j}$ is the coupling strength 
between the $i$-th and $j$-th oscillators. 
The matrix $K = [K_{i,j}]$ may also be viewed as the adjacency matrix for 
the underlying weighted graph.
$\Omega=(\omega_{1},\dots,\omega_{N})$ contains the natural frequencies of the
$N$ oscillators.
Of crucial importance in studying the phase space of this system of ordinary differential equations
are the equilibria which are values of $\theta_1,\dots,\theta_N$ for which
$\frac{d \theta_i}{dt} = 0$ for all $i=1,\dots,N$.
The stability analysis of the equilibria can reveal the behavior
of the dynamical system near the equilibria. 
One can also study synchronization phenomenon of the Kuramoto model with the 
knowledge of the equilibria of the model \cite{KuramotoPaper1}. 
Interpreted as a special case \cite{baillieul1982geometric,dorfler2010synchronization,marecek2014power} of the 
power flow equations in alternating-current power systems
with harmonic currents,
the equilibria of the Kuramoto model provide crucial information for planning and designing power grids.

Equilibria of the Kuramoto model can be found by solving the equilibrium conditions 
$\frac{d\theta_{i}}{d t}=0$ for all $i$.  
Since the equilibrium conditions are invariant under 
transformations of the form $\theta_i \rightarrow \theta_i + \alpha$ for all 
$\theta_{i}$ with any fixed $\alpha\in(-\pi,\pi]$, 
they possess infinitely many solutions. 
To remove this degree of freedom, we fix $\theta_{N}=0$ and remove the $N$-th 
equation $\frac{d \theta_{N}}{dt}=0$.
Thus, we are left with $N-1$ nonlinear equations in $N-1$ angles. 
For example with $N=3$, the Kuramoto model is thus known \cite{tavora1972equilibrium,baillieul1984critical,klos1991physical} to have at most $6$ equilibria on the complete graph.

The counting the number of equilibria of the 
Kuramoto model with a finite number $N$ of oscillators, 
has a very long history and should be seen as a key structural property of the set of equilibria.
Generally, it is being turned into a 
\emph{root counting} problem for systems of polynomial equations.
As we will show, there are multiple reformulations to
polynomial equations, and some of the bounds on the number
of roots are sensitive to the choice of the reformulation.
We present strong numerical evidence that for many Kuramoto models,
a bound 
based on the reformulation we introduce 
and the theory of
\cite{Bernstein75}, \cite{Kushnirenko76}, and \cite{Khovanski78} (BKK)
is the best known and in some sense the best possible bound for the number of complex equilibria.

Outside of the structural properties of the set of equilibria,
one may be interested in the optimisation over the set.
This has been particularly popular in power systems applications,
where there is a cost (of power generation) associated with the equilibria, cf. \cite{Ghaddar2016}.
Again, we demonstrate that the optimisation methods are 
very sensitive to the choice of the reformulation.
Again, the use of the reformulation we introduce offers
performance superior to those previously used.

First, in \S~\ref{Sec:Formulations}, we describe a novel polynomial reformulation of the 
equilibrium equations of the Kuramoto model and compare it the reformulations suggested  previously.
Next, we review some known structural results in \S~\ref{Sec:Existing} abd propose an upper bound on the number of equilibria based on 
our reformulation and the BKK theory in \S~\ref{Sec:BKK}.
Next, we suggest the use of the method of moments for optimisation 
over the set of equilibria in \S~\ref{sec:opt}.
Next, we provide computational results for the Kuramoto model on benchmark 
graphs. 
We show how sensitive is the performance of both  the method of moments and a numerical polynomial homotopy continuation (NPHC) method to the choice of the  reformulation. 
Using the NPHC method, we provide strong
numerical evidence to demonstrate that our structural results are tight in certain cases.
That is: we employ the numerical homotopy continuation method and 
demonstrate that for generically chosen natural frequencies
and weights on the graph, our upper bound is always equal to the number of complex solutions. 
In \S~\ref{sec:Discussion-and-conclusion}, we discuss the implications of our results and conclude.



\section{Polynomial Formulations}\label{Sec:Formulations}

In studying nonlinear systems of equations like the equilibrium equations of
\eqref{Eq:KuramotoStandard}, it is a common practice to first transform
them into algebraic equations 
which would allow the use of powerful tools from algebraic geometry.

\emph{Formulation 1.} Previously  \cite{baillieul1982geometric,Mehta:2009,Mehta:2009zv,Mehta:2011xs,KuramotoPaper1}, 
the equilibrium conditions (\ref{Eq:KuramotoStandard}) were 
transformed into a system of polynomial equations by using the identities 
$\sin (\theta_{i} - \theta_{j}) = \sin \theta_{i} \cos \theta_{j} - \sin \theta_{j} \cos \theta_{i}$
and then the substitution $s_i := \sin \theta_{i}$ and $c_i := \cos \theta_{i}$
for all $i = 1, \dots, N -1$, and adding the 
 equation  
$s_{i}^2 + c_{i}^2 - 1 = 0,$ for all $i = 1, \dots, N-1$. This way, one obtains the following system
of polynomial equations:
\begin{equation}
\tag{F1}
    \label{eq:sincos}
    \begin{aligned}
        \omega_i - \frac{1}{N} \sum_{j=1}^{N} 
            K_{i,j}\left( s_i c_j - s_j c_i \right) &= 0 \\
            s_i^2 + c_i^2 - 1 &= 0,
    \end{aligned}
\end{equation}
for $i = 1, \dots, N-1$.

\emph{Formulation 2.} Alternatively, one may consider a formulation based on the so-called tangent half-angle identities  \cite[pages 382--383]{spivak2006calculus}:
\begin{align}
\sin\theta_{i} = \frac{2 \tan\frac{\theta_{i}}{2}}{1 + \tan^2\frac{\theta_{i}}{2}} \textrm{ and } \cos\theta_{i} = \frac{1- \tan^\frac{\theta_{i}}{2}}{1 + \tan^2\frac{\theta_{i}}{2}},
\end{align}
wherein one introduces  $t_{i}:=\tan \frac{\theta_{i}}{2}$. 
This way, one obtains the following system
of polynomial equations:
\begin{equation}
\tag{F2}
    \begin{aligned}
        \omega_i - \frac{1}{N} \sum_{j=1}^{N} 
            K_{i,j}\left( \frac{2 t_i}{1 + t_i^2} \frac{1 - t_j}{1 + t_j^2} - \frac{2 t_j}{1 + t_j^2} \frac{1 - t_i}{1 + t_i^2} \right) &= 0 \\
            \left( \frac{2 t_i}{1 + t_i^2} \right)^2 + \left( \frac{1 - t_i}{1 + t_i^2} \right)^2 - 1 &= 0,
    \end{aligned}
\end{equation}
for $i = 1, \dots, N-1$.
In Figure \ref{figuretan3}, we present an example of the first equation for $N=3$. 
It is clear that the degree of each equation in this algebrization increases,
as one clears out the denominators to get the equations into the polynomial form,
but 
the number of equations and the number of variables
remains linear in $N$. 

This increase of degree can be partially mitigated by lifting the problem, i.e., 
introducing additional variables such as $a_{i}:= 1 + t_i^2$
and performing the appropriate substitutions.
Thereby, one obtains yet different algebrizations.

We are not aware of any applications of such half-angle formulation to the Kuramoto model, although it has been used under a variety of names in a variety of applications.
It is, for example, known as the Weierstrass substitution, and traceable to the work of Euler \cite[E342, Caput V, paragraph 261]{euler1768institutionum}. 
Emiris  
pioneered the use of this transformation
in Robotics \cite{emiris1994sparse} and Computational Chemistry \cite{emiris1999computer}.

\begin{figure*}[htb]
\begin{equation}
\label{equ:tan3}
\begin{aligned}
  K_{1,2} \left( \frac{2 t_1}{1 + t_1^2} \frac{1 - t_2}{1 + t_2^2} - \frac{2 t_2}{1 + t_2^2} \frac{1 - t_1}{1 + t_1^2} \right) 
+ K_{1,3} \left( \frac{2 t_1}{1 + t_1^2} \frac{1 - t_3}{1 + t_3^2} - \frac{2 t_3}{1 + t_3^2} \frac{1 - t_1}{1 + t_1^2} \right) 
 &= 3 \omega_1  \\
   K_{2,1} \left( \frac{2 t_2}{1 + t_2^2} \frac{1 - t_1}{1 + t_1^2} - \frac{2 t_1}{1 + t_1^2} \frac{1 - t_2}{1 + t_2^2} \right) 
+ K_{2,3} \left( \frac{2 t_2}{1 + t_2^2} \frac{1 - t_3}{1 + t_3^2} - \frac{2 t_3}{1 + t_3^2} \frac{1 - t_2}{1 + t_2^2} \right) 
 &= 3 \omega_1  \\
\end{aligned}
\end{equation}
\caption{A reformulation of the example using $t_{i}:=\tan \frac{\theta_{i}}{2}$.}
\label{figuretan3}
\end{figure*}



\emph{Formulation 3.}
In the present work, we suggest yet another transformation 
using the trigonometric identity \cite[cf. p. 71, 4.3.1]{abramowitz1965}:
\begin{align}
\sin(\theta_{i} - \theta_{j}) = \tfrac{1}{2\mathbb{I}}(e^{\mathbb{I}(\theta_{i} - \theta_{j})} - e^{-\mathbb{I} (\theta_{i} - \theta_{j})})
\end{align}
where $\mathbb{I} := \sqrt{-1}$ is the imaginary unit. The equilibrium
equations of \eqref{Eq:KuramotoStandard} become 
\begin{equation}
	\label{equ:exp}
	\omega_{i} - 
	\sum_{j=1}^{N} 
	\frac{K_{i,j}}{\mathbb{I} N} (
		e^{\mathbb{I} \theta_{i}} e^{-\mathbb{I} \theta_{j}} - 
		e^{-\mathbb{I} \theta_{i}} e^{\mathbb{I} \theta_{j}}
	) = 0.
\end{equation}
To formulate the equilibrium equations as an algebraic system, we let 
\begin{equation} 
	x_{i} := e^{\mathbb{I} \theta_{i}}
	\quad\quad\text{ and }\quad\quad
	y_{i} := e^{-\mathbb{I} \theta_{i}},
	\label{eq:newform}
\end{equation} 
for all $i = 1, \dots, N-1$. 
With this substitution,  \eqref{equ:exp} becomes
%
an enlarged system of $2(N-1)$ equations in $2(N-1)$ variables
\begin{equation}
\tag{F3}
	\begin{aligned}
	\sum_{j=1}^{N} \frac{K_{i,j}}{\mathbb{I} N} 
	(x_i y_j - x_j y_i) &= \omega_{i}
	& \text{ for } i&=1,\dots,N-1 \\
	x_i y_i &= 1 
	& \text{ for } i&=1,\dots,N-1.
	\end{aligned}
	\label{equ:kuramoto-poly}
\end{equation}    
For example, with $N=3$, the system
 \eqref{Eq:KuramotoStandard} 
becomes
\begin{equation}
\label{equ:kuramoto-poly3}
\begin{aligned}
\frac{K_{1,2}}{3\mathbb{I} } (x_1 y_2 - x_2 y_1) + \frac{K_{1,3}}{3\mathbb{I} } (x_1 - y_1)
- \omega_{1} &= 0 \\
\frac{K_{2,1}}{3\mathbb{I}} (x_2 y_1 - x_1 y_2) + \frac{K_{2,3}}{3\mathbb{I}}(x_2 - y_2)
- \omega_{2} &= 0 \\
x_1 y_1 - 1 &= 0 \\
x_2 y_2 - 1 &= 0.
\end{aligned}
\end{equation}
It can be readily verified that the equilibria of \eqref{Eq:KuramotoStandard}
(with the translation symmetry removed) are in one-to-one correspondence
with the special solutions of the above system \eqref{equ:kuramoto-poly} 
that satisfy the additional restriction that $|x_i| = |y_i| = 1$ for 
$i=1,\dots,N-1$.
%
%

It is not clear, however, which of these formulations to use in which applications.
In this paper, 
we consider two criteria related to two applications. 

\section{Bounds on the Number of Equilibria} \label{Sec:Existing}

Via the transformations given in \eqref{eq:sincos} or \eqref{equ:kuramoto-poly}, 
the problem of counting equilibria of the Kuramoto model is turned into a 
\emph{root counting} problem for systems of polynomial equations.
This approach has a long history, going back to 
 \cite{tavora1972equilibrium,baillieul1984critical,klos1991physica}, who have shown that for $N=3$, the Kuramoto model on a complete graph of 3 nodes
has at most $6$ equilibria.
Perhaps the best-known bound on the number of roots of a system of polynomial equations 
comes from the theorem of B\'ezout.
 
\emph{B\'ezout's bound} is simply the product of the degrees of all the 
equations.
In the example shown in \eqref{equ:kuramoto-poly3}, since each of the four 
equations is quadratic (degree 2), the highest possible number of isolated 
solutions as given by B\'ezout's bound is therefore
$2^4 = 16$. 
In general, B\'ezout's bound for the system (\ref{equ:kuramoto-poly}) is
$2^{2(N-1)}$.
B\'ezout's bound is a basic result in \emph{intersection theory} 
\cite{fulton_intersection_1998}, 
the study of how the varieties defined by an ideal given by an algebraic equation intersect one another.

\emph{Bi-homogeneous B\'ezout's bound}
for the
equilibrium equations for \eqref{Eq:KuramotoStandard}
can be derived \cite{baillieul1982geometric,li1987numerical,marecek2014power} as 
$\binom{2(N-1)}{N-1}$.
In a wide variety of special cases of the Kuramoto model, 
there are case-specific bounds as well.
These include  
complete graphs \cite{Casetti:June2003:0022-4715:1091}, 
nearest-neighbour coupling on one-dimensional 
lattice graphs \cite[e.g.]{Casetti:June2003:0022-4715:1091,Mehta:2009,Mehta:2010pe},
two-dimensional lattice graphs \cite{casetti2011microcanonical,nardini2012density,Nerattini:2012pi,Mehta:2009,Mehta:2009zv,Hughes:2012hg},
and three-dimensional lattice graphs \cite{Hughes:2012hg,Mehta:2013iea}, 
and the homogeneous frequencies case \cite{baillieul1982geometric}.
       


Notice that these bounds relate only to isolated roots. Casetti et al.  \cite{Casetti:June2003:0022-4715:1091} have shown that even after fixing 
the trivial zero mode by setting $\theta_{N}=0$, there may exist infinitely many equilibria,
 known as  {\em incoherent manifolds} \cite{strogatz1991stability}. 
Surprisingly, 
such infinite families of equilibria were also shown to exist in the one-, two- and three-dimensional lattice
model with nearest neighbour interaction in ~\cite{Mehta:2009,Mehta2011,Nerattini:2012pi}, 
where it was demonstrated that the number of infinite families
of equilibria grows exponentially in $N$.

\section{The BKK Bound} \label{Sec:BKK}

In the present contribution, we establish an upper bound on the number of 
equilibria for \eqref{Eq:KuramotoStandard} using a novel polynomial formulation
 and the theory of Bernstein \cite{Bernstein75}, Kushnirenko  \cite{Kushnirenko76}, and Khovanski \cite{Khovanski78},
 to whom we refer to using the acronym BKK.
The advantage of this bound over existing bounds is that it takes into 
consideration of the sparsity of the connections in the underlying network.
This marks a significant leap forward from the recent studies of the Kuramoto model and 
the closely related load flow equations for electric-power networks \cite{marecek2014power,KuramotoPaper1,chen_network_2015} from algebraic view points. 


First, let us briefly review the BKK work  \cite{Bernstein75,Khovanski78,Kushnirenko76} on the number of 
isolated non-zero complex solutions which is a refinement of the B\'ezout bound 
that takes into consideration the monomials that appear in the polynomial system:
Given a polynomial, each of its terms give rise to an \emph{exponent vector}.
For instance, for the term $x^3 y^2 z^1$, the exponent vector is simply 
the vector whose entries are the exponents of $x,y$ and $z$, respectively, 
i.e., $(3,2,1)$.
The choice of this ordering is inconsequential as long as it is kept the same 
for each equation.
The set of all exponent vectors derived from the non-zero terms of an polynomial 
equation is called the \emph{support} of that equation. 
For example, if we arrange the variables in the order of $(x_1,y_1,x_2,y_2)$, 
then the supports of the four equations in \eqref{equ:kuramoto-poly3} are
    \begin{align*}
   &  \{(1,0,0,1),(0,1,1,0),(1,0,0,0),(0,1,0,0),(0,0,0,0)\} \\
   &  \{(0,1,1,0),(1,0,0,1),(0,0,1,0),(0,0,0,1),(0,0,0,0)\} \\
&     \{(1,1,0,0),(0,0,0,0)\} \\
  &   \{(0,0,1,1),(0,0,0,0)\}.
    \end{align*}%

A \emph{convex set} is a set of points in which the line segment connecting 
any pair of points in the set also lie in that set.
The \emph{convex hull} of a set is the minimal convex set containing that set.  
For a polynomial, the convex hull of its support is known as the 
\emph{Newton polytope} of that polynomial.
In the study of convex polytopes, the \emph{mixed volume} of several polytopes
is an important concept.
 which can be considered as a generalization of the
concept of volume into the context of several polytopes.
Given $n$ convex polytopes $Q_1,\dots,Q_n \subset \mathbb{R}^n$ and positive 
real numbers $\lambda_1,\dots,\lambda_n$ Minkowski's Theorem 
states that the $n$-dimensional volume of the 
\emph{Minkowski sum}
$\lambda_1 Q_1 + \cdots + \lambda_n Q_n$, defined as
\[
    \{ \;
        \lambda_1 q_1 + \cdots + \lambda_n q_n \,\mid\,
        q_i \in Q_i \text{ for } i = 1,\dots,n \;
    \}
\]
is a homogeneous polynomial of degree $n$ in the variables 
$\lambda_1,\dots,\lambda_n$.
The coefficient associated with the monomial $\lambda_1 \cdots \lambda_n$ 
in this polynomial is known as the \emph{mixed volume} 
of the polytopes $Q_1,\dots,Q_n$.
In the simplest case, the mixed volume of two line segments on the plane is
precisely the area of the parallelogram spanned by translations of these two
line segments.
With these definitions, one can state:

\begin{thm}[Bernstein \cite{Bernstein75}]
    Given a system of $n$ polynomial equations in $n$ variables, the number of 
    isolated complex solutions for which no variable is zero is bounded above 
    by the mixed volume of the Newton polytopes of the equations.
\end{thm}

Recall that via the change of variables given in \eqref{eq:newform}, each
(real) equilibrium of \eqref{Eq:KuramotoStandard} corresponds to a unique
non-zero complex solution of \eqref{equ:kuramoto-poly}.
The BKK bound given above hence provides an upper bound to the
number of isolated (real) equilibria.

Bernstein \cite{Bernstein75} has also shown that the BKK bound is \emph{generically exact}: 
when the coefficients in the polynomial system are chosen 
at random, \emph{with probability one}, the number of isolated complex 
solutions for which no variable is zero is exactly the BKK bound.
In the polynomial formulation of the Kuramoto model given in 
\eqref{equ:kuramoto-poly}, if certain relations are imposed on the coefficients
(e.g., the coefficients of $x_1y_2$ and $x_2y_1$ in \eqref{equ:kuramoto-poly3}
must be the same) the generic exactness still holds true under a mild additional 
condition.
This result of \cite{Bernstein75}, translated to the language of Kuramoto model is thus:

\begin{thm}
    If there exists a choice of $K_{i,j}$'s and $\omega_i$'s for which
    the number of non-zero complex solutions of \eqref{equ:kuramoto-poly} is
    the BKK bound, then for almost all choices of complex $K_{i,j}$'s and 
    $\omega_i$, the number of non-zero complex solutions of 
    \eqref{equ:kuramoto-poly} will be the BKK bound.
\end{thm}

In other words, among the systems \eqref{equ:kuramoto-poly} for all possible
choices of the $K_{i,j}$'s and $\omega_i$'s, if the BKK bound is attainable
then it must also be generically exact.
In Section \ref{Sec:Results}, we shall compute the BKK bound for the polynomial system
\eqref{equ:kuramoto-poly} induced by a number of graphs.
Then, the attainability and hence the generic exactness of the BKK bound in 
each case is verified by solving the system \eqref{equ:kuramoto-poly} for
some specific chosen set of $K_{i,j}$'s and $\omega_i$'s.


\section{The Optimisation Problems} 
\label{sec:opt}

One can also optimise over the system obtained using either of the three
substitutions, possibly intersected with additional polynomial inequalities. 

Let us denote the algebraic set, which is obtained by intersecting the reformulation of 
\eqref{Eq:KuramotoStandard} with $r$ additional polynomial inequalities 
as $K$. Let us consider a polynomial objective function $f(x), x \in K$ and its global minimum $f^*$ achieved at one or more $x^* \in K$. 
Then, it is easy to see:
\begin{thm}[Asymptotic Convergence]
\label{thm:convergent}
Whenever $K$ is non-empty and there exists an $M>0$ such that $\parallel x \parallel_{\infty} < M$ for all $x\in K$,
there exists a hierarchy of semidefinite programming relaxations [SDP$_r$] their respective duals [SDP$_r$]$^*$ such that the following holds: 
\begin{enumerate}
\item[(a) ]
$\inf [SDP_r] \nearrow f^* \text{ as } r\rightarrow\infty,$
\item[(b) ]
$\sup$ of semidefinite-programming duals of SDP$_r]$ $\nearrow f^* \: \text{as} \: r\rightarrow\infty,$
\item[(c) ] if there exists a unique global minimizer $x^* \in K$, with respect to $f$, then as $r$ tends to infinity the components of the optimal solution of [SDP$_r$] corresponding to the linear terms converge to $x^*$.
\end{enumerate}
\end{thm}

The proofs follow from the seminal work of Lasserre \cite{lasserre2001global}.
Although Theorem \ref{thm:convergent} states just the existence of such 
a hierarchy, there are readily available algorithms \cite{Handbook} for constructing the hierarchy and 
computing an arbitrarily-accurate approximation of [SDP$_r$].
Albeit non-trivial, the algorithms have been implemented
successfully \cite[e.g.]{Waki2008,Ghaddar2016}.
Under slightly stronger assumptions \cite{nie2014optimality}, one can also show finite convergence. 

Notice that the inequalities are often bounds on the phase-difference of adjacent oscillators:
\begin{align}
  | \theta_{i} - \theta_{j} | \le z_{ij}, \label{arsinebound}
\end{align}
for some constant $z_{ij}$. 
In power-systems applications, for instance, such inequalities bound losses, thermal limits, and
allow for a certain realism of the equilibria. 
In the traditional reformulation, the inequality \eqref{arsinebound} becomes
\begin{align}
\left( s_i c_j - s_j c_i \right) & \le \arcsin(z_{ij}), \label{arsinebound2}
\end{align}
where $\arcsin(z_{ij})$ is clearly a constant, which can be
precomputed.
In the second reformulation, the constraint \eqref{arsinebound} becomes a bound on $x_i - y_i$.
In the third reformulation, the constraint \eqref{arsinebound} becomes:
\begin{align}
t_i  & \le \arctan(z_{ij}),
\end{align} 
with care needed to consider the appropriate orthant, but again with a constant right-hand side. 
Still, this leaves the question as to which reformulation to use open.

\section{Computational Results} \label{Sec:Results}


First, we compare the BKK bound for both formulations (\ref{eq:sincos}) and 
(\ref{equ:kuramoto-poly}) on the sparsest connected graphs, known as the path graphs.
In a path graph, the $i$-th node for $1 < i < N-1$ is connected to two of
its neighbors: the $(i-1)$-th node and the $(i+1)$-th node forming a path.
The results are presented in Table \ref{tab:path}.

\begin{table*}
    \centering
    \caption{Comparison among different bounds on the number of equilibria and the actual number of complex solutions for generic parameter-values
    for the Kuramoto model on the path graph.}
    \label{tab:path}
    \vspace{0.05in}
    \begin{tabular}{rrrrrrrrrrrrrrr}
        \toprule
        Nodes & 3 & 4 & 5 & 6 & 7 & 8 & 9 & 10 & 11 & 12 & 13 & 14 & 15 \\ 
        \midrule
        B\'ezout's & 16 & 64 & 256 & 1024 & 4096 & 16384 & 65536 & 262144 & 1048576 & 4194304 & 16777216 & 67108864 & 268435456\\ 
        Bi-h. B\'ezout's  & 6 & 20 & 70 & 252 & 924 & 3432 & 12870 & 48620 & 184756 & 705432 & 2704156 & 10400600 & 40116600  \\ 
        BKK for \eqref{eq:sincos} & 8 & 24 & 80 & 256 & 832 & 2688 & 8704 & 28160 & 91136 & 294912 & 954368 & 3088384 & 9994240  \\ 
        BKK for \eqref{equ:kuramoto-poly} & 4 & 8 & 16 & 32 & 64 & 128 & 256 & 512 & 1024 & 2048 & 4096 & 8192 & 16384  \\ 
        Generic root count & 4 & 8 & 16 & 32 & 64 & 128 & 256 & 512 & 1024 & 2048 & 4096 & 8192 & 16384  \\ 
         \bottomrule
    \end{tabular}
\end{table*}

In particular, the generic root count of Table \ref{tab:path} 
presents the results of experiments using a numerical polynomial homotopy continuation method (NPHC). The NPHC method
guarantees that one will obtain all isolated complex solutions for a system of polynomial equations by following the 
following strategy \cite{Li:2003,SW05}:
to solve a system of polynomial equations, one starts with an upper bound on the number of complex solutions of the system.
Then, another system is created such that the system has exactly the same number of complex solutions as the upper bound, and
it is easy to solve. Finally, each solution of this new system is evolved over a single parameter towards the system to be solved.
In particular, we have used the computational packages 
\emph{HOM4PS-3.0} of \cite{Li:03} and \cite{chen2014hom4ps}, as well as  
\emph{Bertini} of \cite{BHSW06,BHSW13}. 

Second, we compare the performance of optimisation methods suggested in Theorem \ref{thm:convergent}
on formulations (\ref{eq:sincos}) and 
(\ref{equ:kuramoto-poly}) in Table~\ref{tab2:path}.
\eqref{equ:kuramoto-poly} turns out to perform better than \eqref{eq:sincos},
due to the numbers of variables being similar in \eqref{equ:kuramoto-poly} and \eqref{eq:sincos},
while the degrees of some of the monomials in \eqref{equ:kuramoto-poly} are lower than in \eqref{eq:sincos}.

More specifically, Table \ref{tab2:path} presents the dimensions of SDP relaxations obtained using SparsePOP of \cite{Waki2008}
at the first applicable level of the hierarchy of Theorem \ref{thm:convergent} 
for generic parameter-values for the Kuramoto
 model on the path graph.
Note that at $r = 1$, the hierarchies of  \cite{Waki2008} and \cite{lasserre2001global} coincide. 
The $m \times n$ constraint matrix $A$ of the SDP relaxation is described by the product $mn$ (dim. of $A$) and 
the number of non-zero entries therein (nnz. of $A$).
These measures influence the memory requirements of any solver.
Additionally, we list the maximum $n_i$ among $n_i \times n_i$ positive semi-definite blocks (max. PSD block),
which influence the run-time of primal-dual interior-point methods such as SeDuMi  \cite{sturm1999}.


\begin{table*}
    \centering
    \caption{Comparison among the dimensions of SDP instances obtained at the first applicable level of the hierarchy of Theorem \ref{thm:convergent} for generic parameter-values for the Kuramoto model on the path graph.}
    \label{tab2:path}    
\begin{tabular}{rrrrrrrrrrrrrrr}\toprule
            Nodes       &  2  &  3  &  4  &   5  &   6  &   7  &   8  &   9  &  10  &  11  &  12  &  13  &  14  &   15  \\ \midrule 
   Dim. of $A$ for \eqref{eq:sincos}  & 200 & 450 & 800 & 1250 & 1800 & 2450 & 3200 & 4050 & 5000 & 6050 & 7200 & 8450 & 9800 & 11250 \\  
   Nnz. of $A$ for \eqref{eq:sincos}  &  20 &  30 &  40 &  50  &  60  &  70  &  80  &  90  &  100 &  110 &  120 &  130 &  140 &  150  \\  
 Max. PSD block for \eqref{eq:sincos} &  3  &  3  &  3  &   3  &   3  &   3  &   3  &   3  &   3  &   3  &   3  &   3  &   3  &   3   \\  
   Dim. of $A$ for \eqref{equ:kuramoto-poly}  & 200 & 450 & 800 & 1250 & 1800 & 2450 & 3200 & 4050 & 5000 & 6050 & 7200 & 8450 & 9800 & 11250 \\  
   Nnz. of $A$ for \eqref{equ:kuramoto-poly}  &  18 &  27 &  36 &  45  &  54  &  63  &  72  &  81  &  90  &  99  &  108 &  117 &  126 &  135  \\ 
 Max. PSD block for \eqref{equ:kuramoto-poly} &  3  &  3  &  3  &   3  &   3  &   3  &   3  &   3  &   3  &   3  &   3  &   3  &   3  &   3   \\ \bottomrule
\end{tabular}
\end{table*}

\section{Discussion and Conclusion}\label{sec:Discussion-and-conclusion}

In this article, we have reformulated the stationary equations of the Kuramoto model to polynomial equations in three different ways.
One of the reformulations is novel. The so-called half-angle transform has been used across several fields, but has not been applied to the Kuramoto model, as far as we know.
Finally, one of the reformulatins is implicit in much related work on the Kuramoto model. 
All three allow for the use of results from semi-algebraic and algebraic geometry. 

In terms of structural results, 
we have provided a prescription to compute an upper bound on the number of equilibria of the model for a given
graph topology, which we called the BKK bound. 
For the complete graph with arbitrary (and inhomogeneous) coupling strengths and 
natural frequencies,
this bound matches the best previously available upper bound, $\binom{2(N-1)}{N-1}$. 
We have demonstrated, however, for sparser graphs such as path graphs, the BKK bound for a new polynomial 
formulation is significantly lower than bi-homogeneous B\'ezout's bound, 
as well as the BKK bound for the traditional polynomial formulation.

We also provide constructive results.
We show how to use the reformulations with the so-called method of moments, which has been developed in semi-algebraic geometry, and which makes it possible to optimise over the stationary equations and a variety of further equalities and inequalities. 
We also demonstrate the computational trade-offs of using the three reformulations.
This may often be preferable to the use of homotopy continuation methods.
Where homotopy continuation methods are used, the BKK bound can be considered as means of constructing 
a starting system in solving the stationary equations.
The bound also provides a concrete stopping criterion to any stochastic method for solving the non-linear equation. 

Considering that systems with sine and cosine of angles and difference of angles are not unique to 
the Kuramoto model, this can have far-reaching implications.
For instance, we point out a remarkable parallel between upper bounds for the equilibria of the Kuramoto model
and those for the equilibria of the complete power flow equations \cite{mehta2015recent,molzahn2015toward}. 
We should also like to point out that BKK-like results have been shown
for affine spaces as well \cite{Li96thebkk,Mau:94,MauW:96,Huber1997},
which could perhaps be used. 




\bibliographystyle{ieeetr}
\bibliography{kuramoto_references,other}

\begin{thebibliography}{10}

\bibitem{kuramoto1975self}
Y.~Kuramoto, ``Self-entrainment of a population of coupled non-linear
  oscillators,'' in {\em International symposium on mathematical problems in
  theoretical physics}, pp.~420--422, Springer, 1975.

\bibitem{acebron2005kuramoto}
J.~A. Acebr{\'o}n, L.~L. Bonilla, C.~J.~P. Vicente, F.~Ritort, and R.~Spigler,
  ``The {K}uramoto model: A simple paradigm for synchronization phenomena,''
  {\em Rev. Mod. Phys.}, vol.~77, no.~1, p.~137, 2005.

\bibitem{strogatz2000kuramoto}
S.~H. Strogatz, ``From {K}uramoto to {C}rawford: exploring the onset of
  synchronization in populations of coupled oscillators,'' {\em Physica D:
  Nonlinear Phenomena}, vol.~143, no.~1, pp.~1--20, 2000.

\bibitem{FD-FB:13b}
F.~D{\"o}rfler and F.~Bullo, ``Synchronization in complex oscillator networks:
  A survey,'' {\em Automatica}, vol.~50, pp.~1539--1564, June 2014.

\bibitem{KuramotoPaper1}
D.~Mehta, N.~S. Daleo, F.~D{\"o}rfler, and J.~D. Hauenstein, ``Algebraic
  geometrization of the {K}uramoto model: Equilibria and stability analysis,''
  {\em Chaos}, vol.~25, no.~5, p.~053103, 2015.

\bibitem{baillieul1982geometric}
J.~Baillieul and C.~I. Byrnes, ``Geometric critical point analysis of lossless
  power system models,'' {\em Circuits and Systems, IEEE Transactions on},
  vol.~29, no.~11, pp.~724--737, 1982.

\bibitem{dorfler2010synchronization}
F.~D{\"o}rfler and F.~Bullo, ``Synchronization and transient stability in power
  networks and non-uniform {K}uramoto oscillators,'' {\em SIAM Journal on
  Control and Optimization}, vol.~50, no.~3, pp.~1616--1642, 2010.

\bibitem{marecek2014power}
J.~Marecek, T.~McCoy, and M.~Mevissen, ``Power flow as an algebraic system,''
  {\em arXiv preprint arXiv:1412.8054}, 2014.

\bibitem{tavora1972equilibrium}
C.~J. Tavora and O.~J.~M. Smith, ``Equilibrium analysis of power systems,''
  {\em IEEE Transactions on Power Apparatus and Systems}, vol.~PAS-91,
  pp.~1131--1137, May 1972.

\bibitem{baillieul1984critical}
J.~Baillieul, ``The critical point analysis of electric power systems,'' in
  {\em The 23rd IEEE Conference on Decision and Control}, no.~23, pp.~154--159,
  1984.

\bibitem{klos1991physical}
A.~Klos and J.~Wojcicka, ``Physical aspects of the nonuniqueness of load flow
  solutions,'' {\em International Journal of Electrical Power \& Energy
  Systems}, vol.~13, no.~5, pp.~268--276, 1991.

\bibitem{Bernstein75}
D.~N. Bernstein, ``The number of roots of a system of equations,'' {\em Funkts.
  Anal. Pril.}, vol.~9, pp.~1--4, 1975.

\bibitem{Kushnirenko76}
A.~G. Kushnirenko, ``Newton polytopes and the {B}ezout theorem,'' {\em Funkts.
  Anal. Pril.}, vol.~10, no.~3, 1976.

\bibitem{Khovanski78}
A.~G. Khovanskii, ``Newton polyhedra and the genus of complete intersections,''
  {\em Funkts. Anal. Pril.}, vol.~12, no.~1, pp.~51--61, 1978.

\bibitem{Ghaddar2016}
B.~Ghaddar, J.~Marecek, and M.~Mevissen, ``Optimal power flow as a polynomial
  optimization problem,'' {\em IEEE T Power. Syst.}, vol.~31, pp.~539--546, Jan
  2016.

\bibitem{Mehta:2009}
D.~Mehta, ``{Lattice vs. Continuum: Landau Gauge Fixing and 't Hooft-Polyakov
  Monopoles},'' {\em Ph.D. Thesis, The Uni. of Adelaide, Australasian Digital
  Theses Program}, 2009.

\bibitem{Mehta:2009zv}
D.~Mehta, A.~Sternbeck, L.~von Smekal, and A.~G. Williams, ``{Lattice Landau
  Gauge and Algebraic Geometry},'' {\em PoS}, vol.~QCD-TNT09, p.~025, 2009.

\bibitem{Mehta:2011xs}
D.~Mehta, ``{Finding All the Stationary Points of a Potential Energy Landscape
  via Numerical Polynomial Homotopy Continuation Method},'' {\em Phys.Rev.},
  vol.~E84, p.~025702, 2011.

\bibitem{spivak2006calculus}
M.~Spivak, ``Calculus,'' {\em Cambridge University Press}, vol.~48, 2006.

\bibitem{euler1768institutionum}
L.~Euler, {\em Calculi integralis}.
\newblock 1768.

\bibitem{emiris1994sparse}
I.~Z. Emiris, {\em Sparse elimination and applications in kinematics}.
\newblock PhD thesis, University of California, Berkeley, 1994.

\bibitem{emiris1999computer}
I.~Z. Emiris and B.~Mourrain, ``Computer algebra methods for studying and
  computing molecular conformations,'' {\em Algorithmica}, vol.~25, no.~2-3,
  pp.~372--402, 1999.

\bibitem{abramowitz1965}
M.~Abramowitz and I.~A. Stegun, {\em Handbook of mathematical functions: with
  formulas, graphs, and mathematical tables}.
\newblock Courier Corporation, 1965.

\bibitem{fulton_intersection_1998}
W.~Fulton, {\em Intersection {{Theory}}}.
\newblock {Springer New York}, Jan. 1998.

\bibitem{li1987numerical}
T.-Y. Li, T.~Sauer, and J.~A. Yorke, ``Numerical solution of a class of
  deficient polynomial systems,'' {\em SIAM journal on numerical analysis},
  vol.~24, no.~2, pp.~435--451, 1987.

\bibitem{Casetti:June2003:0022-4715:1091}
L.~Casetti, M.~Pettini, and E.~G.~D. Cohen, ``Phase transitions and topology
  changes in configuration space,'' {\em Journal of Statistical Physics},
  vol.~111, pp.~1091--1123(33), June 2003.

\bibitem{Mehta:2010pe}
D.~Mehta and M.~Kastner, ``{Stationary point analysis of the one-dimensional
  lattice Landau gauge fixing functional, aka random phase XY Hamiltonian},''
  {\em Annals Phys.}, vol.~326, pp.~1425--1440, 2011.

\bibitem{casetti2011microcanonical}
L.~Casetti, C.~Nardini, and R.~Nerattini, ``Microcanonical relation between
  continuous and discrete spin models,'' {\em Physical Review Letters},
  vol.~106, no.~5, p.~057208, 2011.

\bibitem{nardini2012density}
C.~Nardini, R.~Nerattini, and L.~Casetti, ``Density of states of continuous and
  discrete spin models: a case study,'' {\em Journal of Statistical Mechanics:
  Theory and Experiment}, vol.~2012, no.~02, p.~P02007, 2012.

\bibitem{Nerattini:2012pi}
R.~Nerattini, M.~Kastner, D.~Mehta, and L.~Casetti, ``{Exploring the energy
  landscape of XY models},'' {\em Phys.Rev.}, vol.~E87, no.~3, p.~032140, 2013.

\bibitem{Hughes:2012hg}
C.~Hughes, D.~Mehta, and J.-I. Skullerud, ``{Enumerating Gribov copies on the
  lattice},'' {\em Annals Phys.}, vol.~331, pp.~188--215, 2013.

\bibitem{Mehta:2013iea}
D.~Mehta, C.~Hughes, M.~Schr\"ock, and D.~J. Wales, ``{Potential Energy
  Landscapes for the 2D XY Model: Minima, Transition States and Pathways},''
  {\em J. Chem. Phys.}, vol.~139, p.~194503, 2013.

\bibitem{strogatz1991stability}
S.~H. Strogatz and R.~E. Mirollo, ``Stability of incoherence in a population of
  coupled oscillators,'' {\em Journal of Statistical Physics}, vol.~63,
  no.~3-4, pp.~613--635, 1991.

\bibitem{Mehta2011}
D.~Mehta and M.~Kastner, ``Stationary point analysis of the one-dimensional
  lattice landau gauge fixing functional, aka random phase xy hamiltonian,''
  {\em Annals of Physics}, vol.~In Press, pp.~--, 2011.

\bibitem{chen_network_2015}
T.~{Chen} and D.~{Mehta}, ``On the network topology dependent solution count of
  the algebraic load flow equations,'' {\em IEEE T Power. Syst.}, vol.~33,
  pp.~1451--1460, March 2018.

\bibitem{lasserre2001global}
J.~B. Lasserre, ``Global optimization with polynomials and the problem of
  moments,'' {\em SIAM Journal on Optimization}, vol.~11, no.~3, pp.~796--817,
  2001.

\bibitem{Handbook}
M.~F. Anjos and J.-B. Lasserre, eds., {\em Handbook on semidefinite, conic and
  polynomial optimization}, vol.~166 of {\em International series in operations
  research \& management science}.
\newblock New York: Springer, 2012.

\bibitem{Waki2008}
H.~Waki, S.~Kim, M.~Kojima, M.~Muramatsu, and H.~Sugimoto, ``Algorithm 883:
  Sparsepop---a sparse semidefinite programming relaxation of polynomial
  optimization problems,'' {\em ACM Trans. Math. Softw.}, vol.~35,
  pp.~15:1--15:13, July 2008.

\bibitem{nie2014optimality}
J.~Nie, ``Optimality conditions and finite convergence of {L}asserre's
  hierarchy,'' {\em Mathematical programming}, vol.~146, no.~1-2, pp.~97--121,
  2014.

\bibitem{Li:2003}
T.~Y. Li, ``Solving polynomial systems by the homotopy continuation method,''
  {\em Handbook of numerical analysis}, vol.~XI, pp.~209--304, 2003.

\bibitem{SW05}
A.~Sommese and C.~Wampler, {\em The Numerical Solution of Systems of
  Polynomials Arising in Engineering and Science}.
\newblock World Scientific Publishing, Hackensack, NJ, 2005.

\bibitem{Li:03}
T.~L. Lee, T.~Y. Li, and C.~H. Tsai, ``Hom4ps-2.0, a software package for
  solving polynomial systems by the polyhedral homotopy continuation method,''
  {\em Computing}, vol.~83, pp.~109--133, 2008.

\bibitem{chen2014hom4ps}
T.~Chen, T.-L. Lee, and T.-Y. Li, ``Hom4ps-3: A parallel numerical solver for
  systems of polynomial equations based on polyhedral homotopy continuation
  methods,'' in {\em Mathematical Software--ICMS 2014}, pp.~183--190, Springer,
  2014.

\bibitem{BHSW06}
D.~J. Bates, J.~D. Hauenstein, A.~J. Sommese, and C.~W. Wampler. Software
  available at www.nd.edu/$\sim$sommese/bertini, 2006.

\bibitem{BHSW13}
D.~Bates, J.~Hauenstein, A.~Sommese, and C.~Wampler, {\em Numerically solving
  polynomial systems with {B}ertini}, vol.~25.
\newblock SIAM, 2013.

\bibitem{sturm1999}
J.~F. Sturm, ``Using sedumi 1.02, a matlab toolbox for optimization over
  symmetric cones,'' {\em Optimization methods and software}, vol.~11, no.~1-4,
  pp.~625--653, 1999.

\bibitem{mehta2015recent}
D.~Mehta, D.~K. Molzahn, and K.~Turitsyn, ``Recent advances in computational
  methods for the power flow equations,'' in {\em 2016 American Control
  Conference (ACC)}, pp.~1753--1765, July 2016.

\bibitem{molzahn2015toward}
D.~K. Molzahn, D.~Mehta, and M.~Niemerg, ``Toward topologically based upper
  bounds on the number of power flow solutions,'' in {\em 2016 American Control
  Conference (ACC)}, pp.~5927--5932, July 2016.

\bibitem{Li96thebkk}
T.~Y. Li and X.~Wang, ``The {BKK} root count in {$C^{n}$},'' {\em Math. Comp},
  vol.~65, pp.~1477--1484, 1996.

\bibitem{Mau:94}
J.~M. Rojas, ``A convex geometric approach to counting the roots of a
  polynomial system,'' {\em Theor. Comput. Sci.}, vol.~133, no.~1,
  pp.~105--140, 1994.

\bibitem{MauW:96}
J.~M. Rojas and X.~Wang, ``Counting affine roots of polynomial systems via
  pointed newton polytopes,'' {\em J. Complex.}, vol.~12, no.~2, pp.~116--133,
  1996.

\bibitem{Huber1997}
B.~Huber and B.~Sturmfels, ``Bernstein's theorem in affine space,'' {\em
  Discrete Comput. Geom.}, vol.~17, no.~2, pp.~137--141, 1997.

\end{thebibliography}

\end{document}